\documentclass[reqno]{amsart}

\usepackage{pstricks,amsmath,amsthm,amssymb,amsfonts}
\numberwithin{equation}{section} 
\theoremstyle{plain}
\newtheorem{theo+}           {Theorem}      [section]
\newtheorem{prop+}  [theo+]  {Proposition}
\newtheorem{coro+}  [theo+]  {Corollary}
\newtheorem{lemm+}  [theo+]  {Lemma}
\newtheorem{defi+}  [theo+]  {Definition}

\theoremstyle{definition}
\newtheorem{rema+}  [theo+]  {Remark}
\newtheorem{prob+}  [theo+]  {Problem}
\newtheorem{exam+}  [theo+]  {Example}

\newenvironment{theorem}{\begin{theo+}}{\end{theo+}}
\newenvironment{proposition}{\begin{prop+}}{\end{prop+}}
\newenvironment{corollary}{\begin{coro+}}{\end{coro+}}
\newenvironment{lemma}{\begin{lemm+}}{\end{lemm+}}

\newenvironment{remark}{\begin{rema+}}{\end{rema+}}
\newenvironment{definition}{\begin{defi+}}{\end{defi+}}

\newcommand{\ga}{\gamma}

\newcommand{\la}{\lambda}

\newcommand{\om}{\omega}

\newcommand{\tha}{\theta}

\newcommand{\End}{\operatorname{End}}
\newcommand{\ot}{\otimes}

\begin{document}

\baselineskip 18pt
\larger[2]
\title
[Izergin--Korepin identity for 8VSOS model]
{An  Izergin--Korepin-type identity\\ for the 8VSOS model,
with applications to\\ alternating sign matrices
} 
\author{Hjalmar Rosengren}
\address
{Department of Mathematical Sciences
\\ Chalmers University of Technology\\SE-412~96 G\"oteborg, Sweden}
\address
{Department of Mathematical Sciences
\\ 
 University of Gothenburg\\SE-412~96 G\"oteborg, Sweden}

\email{hjalmar@math.chalmers.se}
\urladdr{http://www.math.chalmers.se/{\textasciitilde}hjalmar}

\thanks{Research  supported by the Swedish Science Research
Council (Vetenskapsr\aa det)}

\begin{abstract} We obtain a new expression for the partition function of the 8VSOS model with domain wall boundary conditions, which we 
 consider  to be the natural extension of the Izergin--Korepin formula for the six-vertex model. As applications, we find dynamical (in the sense of the dynamical Yang--Baxter equation) generalizations of the  enumeration and $2$-enumeration of alternating sign matrices. The dynamical enumeration has a nice interpretation in terms of three-colourings of the square lattice. 
\end{abstract}

\maketitle        

\section{Introduction}

An alternating sign matrix is a square matrix with entries $0$, $-1$ and $1$, such that the non-zero entries in each row and column form an alternating sequence of the form
 $$1,-1,1,-1,\dots,-1,1. $$
Mills, Robbins and Rumsey \cite{mrr} conjectured that the number of $n\times n$ alternating sign matrices equals
$$A_n=\frac{1!\, 4!\, 7!\dotsm (3n-2)!}{n!(n+1)!\dotsm (2n-1)!}. $$
This was proved thirteen years later by Zeilberger \cite{z}.
Kuperberg \cite{k} found a simpler proof based on the six-vertex model on a square with domain wall boundary conditions. This is a lattice model of statistical mechanics, whose states can be identified with alternating sign matrices. 
Each  state carries a weight, in general depending on $2n+1$ parameters $q,x_1,\dots,x_n,y_1,\dots,y_n$.  The partition function for the model is the
sum of the weight of all states. By 
the  Izergin--Korepin identity \cite{i,ik}, it can be expressed  in terms of the determinant  
$$\det_{1\leq i,j\leq n}\left(\frac{1}{(x_i-qy_j)(x_i-q^{-1}y_j)}\right).$$

Kuperberg observed that 
  when $q=e^{2\pi i/3}$ and $x_i=y_i=1$ for all $i$, 
 the weight of each state can be normalized to $1$, so the partition function is  equal to $A_n$.  The Izergin--Korepin identity then gives
$$A_n=3^{\binom{n+1}{2}}\lim_{\substack{x_1,\dots,x_n\rightarrow 1\\ y_1,\dots,y_n\rightarrow 1 }}
\prod_{1\leq i<j\leq n}\frac 1{(x_i-x_j)(y_i-y_j)}\det_{1\leq i,j\leq n}\left(\frac{x_i-y_j}{x_i^3-y_j^3}\right).$$
Although computing the limit is  not trivial, Kuperberg could do it by elementary means.

The eight-vertex model is a generalization of the six-vertex model,
where the weights may be taken as elliptic functions of the parameters. In his solution of the eight-vertex model, Baxter \cite{b} introduced a
different generalization of the six-vertex model, 
the 8VSOS (eight-vertex-solid-on-solid)
model. Actually, the  ``8'' is somewhat misleading, since the
model obeys the  ice rule and thus admits
 only six  local configurations.
In particular, imposing domain wall
 boundary conditions,  states can  be identified 
with alternating sign matrices. Compared to the six-vertex model, its main distinguishing feature is the presence of a ``dynamical'' parameter,  associated to the faces of the underlying lattice. 

The purpose of the present paper is to study generalizations of the Izergin--Korepin identity and of Kuperberg's specialization, when the six-vertex model is replaced by the 8VSOS model. One motivation is to understand
the significance of the dynamical parameter from a combinatorial viewpoint. We also hope that our results may be useful for studying the thermodynamic limit of the partition function, similarly as the  Izergin--Korepin identity is used in \cite{kz,zj}.

The plan of the paper is as follows. \S \ref{tfs} contains preliminaries on theta functions, and in \S \ref{ms} we recall the 8VSOS model and obtain some elementary properties of its partition function.
In \S \ref{wfs} we show that the partition function can be identified with a special case of the elliptic weight functions of Tarasov and Varchenko \cite{tv}.
We then give our main result, Theorem \ref{mt}, which  expresses the partition function as a sum of $2^n$ determinants. We argue that this is a natural extension of the Izergin--Korepin identity. The rest of the paper is concerned with the case when $q$ is a root of unity. 
In \S \ref{sds} we show that if $q^N=1$, our generalized Izergin--Korepin identity can be reduced to a sum of $N-1$ determinants.
After the preliminary 
 \S \ref{kss}, we consider the analogue of Kuperberg's specialization 
for the  trigonometric 8VSOS model in  \S \ref{es}. 
Curiously, the corresponding extension of the alternating sign matrix theorem, Corollary \ref{ec}, contains not only the numbers $A_n$, but also the numbers
$$C_n=\prod_{j=1}^n\frac{(3j-1)(3j-3)!}{(n+j-1)!},$$
which enumerate cyclically symmetric plane partitions in a cube of size $n$ \cite{a}. 
This dynamical enumeration is best understood in terms of three-colourings of the square lattice. In  particular, it allows us to compute exactly the probability that a random square from a random three-colouring, satisfying  domain wall boundary conditions, has any given colour, see Corollary \ref{rcc}.
Finally, in  \S \ref{tes} we study the more elementary case of  $2$-enumeration.

During the completion of our work, several closely related papers have appeared.  Pakuliak, Rubtsov and Silantyev \cite{prs} studied the 8VSOS partition function using algebraic techniques, and independently obtained Theorem \ref{wft}. 
 Foda,  Wheeler and Zuparic \cite{fwz} studied some other elliptic models, where the partition function can be explicitly factored. Proposition \ref{sqp} can also be obtained from their results. Finally, 
 Razumov and Stroganov \cite{ras} studied the partition function of the three-colouring model using an elliptic function parametrization. 
Trying to understand the relation to their work greatly improved the exposition in  \S \ref{es}. 

{\bf Acknowledgement:} I would like to thank Professors Michio Jimbo,
 Alexander Razumov, Yuri Stroganov, Vitaly Tarasov and Ole Warnaar 
for vital comments at various stages of the  work.

\section{Theta functions}\label{tfs}
Throughout, $\tau$ and $\eta$ will be fixed parameters such that
 $\operatorname{Im}(\tau)>0$, $\eta\notin\mathbb Z+\tau\mathbb Z$.
We will write 
 $p=e^{2\pi i\tau}$ and $q=e^{2\pi i\eta}$. 
By $q^x$ we  always mean $e^{2\pi i\eta x}$.

We will use the notation
$$[x]=q^{-x/2}\prod_{j=0}^\infty(1-p^jq^x)(1-p^{j+1}q^{-x}). $$
Up to a multiplicative constant,  $[x]$ equals the Jacobi theta function
$\theta_1(\eta x|\tau)$ \cite{ww}. We sometimes write for short
$$[x_1,\dots,x_n]=[x_1]\dotsm[x_n]. $$
 
The function
$x\mapsto[x]$ is  odd, entire, and satisfies
\begin{equation}\label{taf}[x+u,x-u,y+v,y-v]-[x+v,x-v,y+u,y-u]=[x+y,x-y,u+v,u-v].\end{equation}
In fact,  up to an elementary multiplier, the only such function is the Jacobi theta function, together with the degenerate cases 
$[x]=\sin(\pi\eta x)$
and $[x]=x$, cf.\ \cite[p.\ 461]{ww}. We find it helpful to think of $[x]$ as a two-parameter deformation of the number $x$.

In the  case  $q^N=1$, we find it more convenient to  use 
the notation
$$\theta(x)=\theta(x;p)= \prod_{j=0}^\infty(1-p^jx)(1-p^{j+1}/x), $$
$$\theta(x_1,\dots,x_n)=\theta(x_1,\dots,x_n;p)=\theta(x_1;p)\dotsm\tha(x_n;p), $$
so that
$$[x]=q^{-x/2}\theta(q^x;p). $$

The following terminology will be useful.

\begin{definition}
Fixing $\tau$ and $\eta$, we say that $f$ is 
a theta function of order $n$ and norm $t$ if there exist constants
$a_1,\dots, a_n$ and $C$ with $a_1+\dots+a_n=t$, such that
\begin{equation}\label{tff}f(x)=C[x-a_1]\dotsm[x-a_n].\end{equation}
Equivalently, $f$ is an entire function such that 
\begin{equation}\label{qp}f(x+1/\eta)=(-1)^nf(x),\qquad f(x+\tau/\eta)=(-1)^ne^{2\pi i\eta(t-nx)-\pi i\tau n }f(x). \end{equation}
\end{definition}

The equivalence of these two properties is  classical  \cite[p.\ 45]{we}.  More generally, any function of the form
$$f(x)=\sum_j\frac{[x-a_1^{(j)}]\dotsm [x-a_{m+n}^{(j)}]}{[x-b_1^{(j)}]\dotsm [x-b_{m}^{(j)}]},  $$
where 
$$a_1^{(j)}+\dots +a_{m+n}^{(j)}-b_1^{(j)}-\dots -b_{m}^{(j)}=t $$
for each $j$, satisfies the quasi-periodicity \eqref{qp}.
If  $f$ is entire (that is, the singularities at $x=b_i^{(j)}$ are all removable), it can then be factored as in \eqref{tff}.
Unless $f$ is  identically zero, the zero set is then
  $a_i+\mathbb Z\eta^{-1}+\mathbb Z\tau\eta^{-1}$, $1\leq i\leq n$, where 
$a_1+\dots+a_n=t$. Thus, to prove that $f$ vanishes identically, it suffices to find $n$ independent zeroes. This gives a powerful (and classical) method for proving theta function identities, which we are going to apply repeatedly.

Finally, we recall Frobenius' determinant evaluation \cite{f}
\begin{multline}\label{fda}\det_{1\leq i,j\leq n}\left(\frac{[x_i-y_j+t]}{[x_i-y_j]}\right)\\
=\frac{(-1)^{\binom n2}[t]^{n-1}\big[|x|-|y|+t\big]\prod_{1\leq i<j\leq n}[x_j-x_i,y_j-y_i]}{\prod_{i,j=1}^n[x_i-y_j]},\end{multline}
or equivalently
\begin{multline}\label{fd}
\det_{1\leq i,j\leq n}\left(\frac{\tha(tx_i/y_j)}{\tha(x_i/y_j)}\right)\\
=\frac{(-1)^{\binom n2}\tha(t)^{n-1}Y\tha(tX/Y)\prod_{1\leq i<j\leq n}x_jy_j\tha(x_i/x_j,y_i/y_j)}{\prod_{i,j=1}^n y_j\tha(x_i/y_j)}.
\end{multline}
Here and throughout, we write 
\begin{equation}\label{nn}|x|=x_1+\dots +x_n,\qquad X=x_1\dotsm x_n. 
\end{equation}

\section{The 8VSOS model}  
\label{ms}

We will study the 8VSOS model on a square with domain wall boundary conditions. 
There are several ways to describe this type of models, see e.g.\ \cite{p}. 
We find it convenient to use \emph{height  matrices}. Fixing a non-negative integer $n$, by a \emph{state} we mean an  $(n+1)\times(n+1)$ matrix,  such that any two horizontally or vertically adjacent entries differ by $1$, and such that the boundary entries are specified as
\begin{equation}\label{dwbc}\left(\begin{matrix} 0 & 1 & 2 &\dotsm & n\\
1 & &&& n-1\\
2 &&&& n-2\\
\vdots&&&&\vdots\\
n & n-1 & n-2 & \dotsm & 0
\end{matrix}\right). \end{equation}
As an example, when $n=2$ there are two states:
\begin{equation}\label{sss}\left(\begin{matrix}0 & 1 & 2\\ 1&0&1\\2&1&0\end{matrix}\right),\qquad 
\left(\begin{matrix}0 & 1 & 2\\ 1&2&1\\2&1&0\end{matrix}\right). \end{equation}

By a \emph{block}, we mean a
 $(2\times 2)$-block  
 of adjacent matrix entries in a state. The blocks can 
be viewed as entries of an  $n\times n$ matrix. Replacing each block $\left(\begin{smallmatrix}a&b\\c&d\end{smallmatrix}\right)$
by $(b+c-a-d)/2$ gives a bijection between states 
and 
 alternating sign matrices of size $n\times n$. 
For instance, the states \eqref{sss} correspond to the  alternating sign matrices
$$\left(\begin{matrix} 1 & 0\\ 0&1\end{matrix}\right),\qquad 
\left(\begin{matrix}0 & 1 \\ 1&0\end{matrix}\right). $$

We will consider Boltzmann weights labelled by 
 $a,b,c,d\in\{\pm 1\}=\{\pm\}$, satisfying the ``ice rule'' $a+b=c+d$. 
For any such labels, let there be given  a meromorphic 
function $R^{ab}_{cd}$ of two complex variables.   
We also fix $2n+1$ generic complex 
parameters $\lambda,x_1,\dots,x_n,y_1,\dots,y_n$.

If  a block $\left(\begin{smallmatrix}a&b\\c&d\end{smallmatrix}\right)$ has coordinates $(i,j)$, chosen with $1\leq i,j\leq n$ in the standard way, that 
block is said to have local weight
$$R^{b-a,d-b}_{d-c,c-a}(\la+a,x_i-y_j). $$
This describes a generalized ice model, where ``generalized'' refers to the ``dynamical'' or ``face'' parameter $\lambda$, which is absent in  the six-vertex model. 
The weight of a state is defined as the product of all  local weights, and the partition function as
\begin{equation}\label{ss}Z_n(x;y;\lambda)=\sum_{\text{states}}\operatorname{weight}(\text{state}).\end{equation}

As an example, from \eqref{sss} we see that
\begin{equation*}\begin{split}    Z_2(x;y;\lambda)&=R^{+-}_{-+}(\la,x_1-y_1)R^{+-}_{+-}(\la+1,x_1-y_2)\\
&\qquad\times R^{-+}_{-+}(\la+1,x_2-y_1)R^{+-}_{-+}(\la,x_2-y_2)\\
&\quad+
R^{++}_{++}(\la,x_1-y_1)R^{+-}_{-+}(\la+1,x_1-y_2)\\
&\qquad\times R^{+-}_{-+}(\la+1,x_2-y_1)R^{--}_{--}(\la+2,x_2-y_2).
\end{split}\end{equation*}

One is particularly interested in Boltzmann weights satisfying the quantum dynamical Yang--Baxter equation (or star-triangle relation), which can be described as follows \cite{fe}. Let $V=Ve_+\oplus Ve_-$ be a two-dimensional complex vector space, and introduce the operators $R(\la,u)\in\End(V\ot V)$ by  
$$R(\la,u)(e_a\ot e_b)=\sum_{c+d=a+b}R^{ab}_{cd}(\la,u)(e_c\ot e_d). $$
Then, 
\begin{multline}\label{qdyb}R^{12}(\la+h^{3},u_1-u_2)R^{13}(\lambda,u_1-u_3)R^{23}(\lambda+h^{1},u_2-u_3)\\
=R^{23}(\lambda,u_2-u_3)R^{13}(\lambda+h^{2},u_1-u_3)R^{12}(\la,u_1-u_2). 
\end{multline}
This should be understood as an identity for meromorphic functions of $\lambda,u_1,u_2,u_3$ with values in $\End(V^{\ot 3})$, with  notation as explained by the example
$$R^{12}(\la+h^{3},u_1-u_2)(e_a\ot e_b\ot e_c)
=R(\la+c,u_1-u_2)(e_a\ot e_b)\ot e_c.
 $$
Part of the interest in this case comes from the following fundamental fact, which follows from the discussion in \cite[\S 9.6]{b2}.

\begin{proposition}[Baxter]\label{sp}
If  the Boltzmann weights satisfy \eqref{qdyb},  then the partition function $Z_n(x;y;\la)$ is a symmetric function of $x$ and $y$. 
\end{proposition}

Adopting the normalization of \cite{djkmo},
the 8VSOS model is given by the following  solution of \eqref{qdyb}:
$$R^{++}_{++}(\la,u)=R^{--}_{--}(\la,u)=\frac{[u+1]}{[1]}, $$
$$R^{+-}_{+-}(\la,u)=\frac{[u][\la+1]}{[1][\la]},\qquad
R^{-+}_{-+}(\la,u)=\frac{[u][\la-1]}{[1][\la]},
 $$
$$R^{-+}_{+-}(\la,u)=\frac{[\la+u]}{[\la]},\qquad R^{+-}_{-+}(\la,u)=\frac{[\la-u]}{[\la]}. $$
From now on, we restrict our attention to this model.

\begin{lemma}\label{tfc}
The partition function is a theta function of each $x_i$ of order $n$ and norm $|y|+\lambda$, and of each $y_i$ of order $n$ and norm $|x|-\lambda$. 
\end{lemma}

\begin{proof} By Proposition \ref{sp}, it suffices to consider the case $i=1$.
Thus, we consider $f$ as a function of $x_1$; the case of $y_1$ is treated similarly. 
 It is well-known and easy to see that 
for each state
there exists a $k$, with $1\leq k\leq n$, such that the second row of the height matrix is 
$$1\ 2\ 3\ \dotsm \ k\ k-1\ k\ k+1\ \dotsm\ n-1. $$
 The $x_1$-dependent part of the partition function is then
\begin{equation}\label{fr}[x_1-y_1+1]\dotsm [x_1-y_{k-1}+1][\la+k-1-x_1+y_k][x_1-y_{k+1}]\dotsm[x_1-y_n], \end{equation}
 which is a theta function of the desired form.
\end{proof}

\begin{lemma}\label{sl}
The partition function satisfies
\begin{align*}
Z_n(x;y;\lambda)\Big|_{x_1+1=y_1}
&=\frac{[\la+n]\prod_{k=2}^n[y_1-y_k-1][x_k-y_1]}{[\la+n-1][1]^{2(n-1)}}\\
&\quad\times Z_{n-1}(x_2,\dots,x_n;y_2,\dots,y_n;\lambda),\\
Z_n(x;y;\lambda)\Big|_{x_1=y_1}
&=\frac{\prod_{k=2}^{n}[y_1-y_k+1][x_k-y_1+1]}{[1]^{2(n-1)}}\\
&\quad\times Z_{n-1}(x_2,\dots,x_{n};y_2,\dots,y_{n};\lambda+1).
 \end{align*}
\end{lemma}

By Proposition \ref{sp}, there are similar identities for any specialization $x_i=y_j$ and $x_i=y_j+1$, $1\leq i,j\leq n$.

\begin{proof} We first consider the case $x_1+1=y_1$. We observe that if  $k\neq 1$, then \eqref{fr} vanishes. Thus,  only states with $k=1$ contribute to the partition function. 
This fixes also the second column, all such  states having the form
\begin{equation*}\left(\begin{matrix} 0 & 1 & 2 &3&\dotsm & n\\
1 &0 &1&2&\dotsm & n-1\\
2 &1&&&& n-2\\
3 &2&&&& n-3\\
\vdots&\vdots &&&&\vdots\\
n & n-1 & n-2 &n-3& \dotsm & 0
\end{matrix}\right). \end{equation*}
It follows that the partition function factors as
\begin{multline*}R^{+-}_{-+}(\la,-1)\prod_{k=2}^n R^{+-}_{+-}(\la+k-1,x_1-y_k)
R_{-+}^{-+}(\la+k-1,x_k-y_1)\\
\times Z_{n-1}(x_2,\dots,x_n;y_2,\dots,y_n;\lambda),
\end{multline*}
which simplifies to the desired expression.

To prove the second identity it is better to let
$x_1=y_n$, which is equivalent in view of  Lemma \ref{sl}. One may then apply  a similar argument as before.
\end{proof}

\section{Relation to elliptic weight functions}
\label{wfs}

In \cite{tv}, Tarasov and Varchenko introduced \emph{elliptic weight functions}, which have played a fundamental role for constructing solutions to the $q$KZ and $q$KZB equations, see further \cite{ftv1,ftv2}.
The following 
result shows that the partition function for the 8VSOS model is an elliptic weight function. Indeed, if we let $\ell =n$ and $\xi_1=\dots=\xi_n=\sqrt{\eta}$ in  \cite[Eq.\ (2.20)]{tv}, it is straight-forward to identify the two expressions.

\begin{theorem}\label{wft}
The partition function can be represented as
\begin{multline}\label{wfi}Z_n(x;y;\la)=\frac{\prod_{i,j=1}^n[y_j-x_i]}{[1]^{n(n-1)}\prod_{j=1}^n[\la+j-1]}\\
\times\sum_{\sigma\in S_n}\prod_{1\leq i<j\leq n}\frac{[y_{\sigma(j)}-y_{\sigma(i)}+1][y_{\sigma(j)}-x_i-1]}{[y_{\sigma(j)}-y_{\sigma(i)}][y_{\sigma(j)}-x_i]}\\
\times\prod_{j=1}^n\frac{[y_{\sigma(j)}-x_j+\la+n-j]}{[y_{\sigma(j)}-x_j]}.
 \end{multline}
\end{theorem}

\begin{proof} Consider the two sides  of \eqref{wfi} as  functions of $x_1$. 
By Lemma~\ref{tfc}, the left-hand side is a theta function of
 order $n$ and norm $|y|+\lambda$, and it is straight-forward to check that the same is true for the right. Thus, as discussed in \S \ref{tfs}, it suffices to verify \eqref{wfi} for  $x_1=y_j-1$, $1\leq j\leq n$. Since both sides are symmetric in $y$, it is in fact enough to take $x_1=y_1-1$. On the left-hand side, we may then apply Lemma~\ref{sl}. On the right, only terms with $\sigma(1)=1$ are non-zero, so it can be viewed as a sum over $S_{n-1}$. In this way,  
\eqref{wfi} is reduced to the same identity with $n$ replaced by $n-1$, and is thus proved  by induction on $n$. 
\end{proof}

The same result was independently  obtained by Pakuliak et al.\ \cite{prs}.

\section{An extension of the Izergin--Korepin identity}
\label{iks}

The following identity is our main result. Originally, we derived it from Theorem~\ref{wft} by a  complicated argument. However, as soon as one has guessed the formula, it is easy to prove directly.

\begin{theorem}\label{mt}
For generic $\gamma$, 
the partition function can be represented as
\begin{multline}\label{mti}Z_n(x;y;\la)=\frac{(-1)^{\binom n2}[\la+n]}{[1]^{n^2}[\ga]^n\big[|x|-|y|+\la+\ga+n\big]}\\
\times\frac{\prod_{i,j=1}^n[x_i-y_j][x_i+1-y_j]}{\prod_{1\leq i<j\leq n}[x_i-x_j][y_i-y_j]}\\
\times\sum_{S\subseteq\{1,\dots,n\}}(-1)^{|S|}
\frac{[\la+\ga+n-|S|]}{[\la+n-|S|]}\det_{1\leq i,j\leq n}\left(\frac{[x_i^S-y_j+\ga]}{[x_i^S-y_j]}\right),\end{multline}
where
\begin{equation}\label{xs}x_i^S=\begin{cases}x_i+1,& i\in S,\\ x_i,& i\notin S.\end{cases} \end{equation}
\end{theorem}

Before we prove  Theorem \ref{mt}, we make some remarks.

\begin{remark}\label{fdr}
The determinants in \eqref{mti} are evaluated by \eqref{fda}. 
This leads to the equivalent identity
\begin{multline}\label{wi}Z_n(x;y;\la)
=\frac{[\la+n]}{[1]^{n^2}[\ga]\big[|x|-|y|+\la+\ga+n\big]}\\
\times\sum_{S}(-1)^{|S|}
\frac{[\la+\ga+n-|S|]\big[|x|-|y|+\ga+|S|\big]}{[\la+n-|S|]}\\
\times\prod_{i\in S,\,j\notin S}\frac{[x_i+1-x_j]}{[x_i-x_j]}\prod_{j=1}^n\left(\prod_{i\in S}[x_i-y_j]\prod_{i\notin S}[x_i+1-y_j]\right),\end{multline}
which expresses $Z_n$ as a sum of $2^n$ explicitly factored terms. This  is much  better than the $A_n$ terms in \eqref{ss} or the $n!$ terms in \eqref{wfi}. In the case of the six-vertex model, the corresponding identity (without the  freedom of choosing $\gamma$) is discussed by Warnaar \cite[Eq.\ (3.3)]{wa}.
\end{remark}

\begin{remark}
Theorem \ref{mt} can be viewed as a generalization of the
Izergin--Korepin identity. To see this, consider the degenerate case when $[x]=\sin(\pi\eta x)$, with $\operatorname{Im}(\eta)>0$, 
and $\lambda\rightarrow \infty$. Up to normalization, this limit corresponds to the six-vertex model.  Since  $[\la+a]/[\la+b]\rightarrow e^{i\pi \eta(b-a)}$, we obtain in the limit
\begin{multline*}Z_n(x;y;\infty)=\frac{(-1)^{\binom n2}e^{i\pi\eta(|x|-|y|)}}{[1]^{n^2}[\ga]^n}
\frac{\prod_{i,j=1}^n[x_i-y_j][x_i+1-y_j]}{\prod_{1\leq i<j\leq n}[x_i-x_j][y_i-y_j]}\\
\times\sum_{S\subseteq \{1,\dots,n\}}(-1)^{|S|}
\det_{1\leq i,j\leq n}\left(\frac{[x_i^S-y_j+\ga]}{[x_i^S-y_j]}\right).
\end{multline*}
By linearity of the determinant and a trigonometric identity, the sum in $S$ can  be written as the single determinant
\begin{multline*}\det_{1\leq i,j\leq n}\left(\frac{[x_i-y_j+\ga]}{[x_i-y_j]}
-\frac{[x_i+1-y_j+\ga]}{[x_i+1-y_j]}
\right)\\
=\det_{1\leq i,j\leq n}\left(\frac{[1][\ga]}{[x_i-y_j][x_i+1-y_j]}
\right).
\end{multline*}
Thus, $\gamma$ cancels and we obtain the Izergin--Korepin identity in the form
\begin{multline*}Z_n(x;y;\infty)=\frac{(-1)^{\binom n2}e^{i\pi\eta(|x|-|y|)}}{[1]^{n^2-n}}
\frac{\prod_{i,j=1}^n[x_i-y_j][x_i+1-y_j]}{\prod_{1\leq i<j\leq n}[x_i-x_j][y_i-y_j]}\\
\times\det_{1\leq i,j\leq n}\left(\frac{1}{[x_i-y_j][x_i+1-y_j]}
\right).
\end{multline*}
\end{remark}

\begin{proof}[Proof of Theorem \ref{mt}] When $n=1$, 
the result can be simplified as
\begin{multline*}[\la+\ga+1][x-y+\ga][\la][x-y+1]-[\la+\ga][x-y+\ga+1][\la+1][x-y]\\
=[1][\ga][y-x+\la][x-y+\la+\ga+1],\end{multline*}
which is an instance of \eqref{taf}.

We now proceed by induction on $n$. 
Let $f_L$ and $f_R$ denote the left-hand and right-hand sides of \eqref{mti} multiplied by  $[|x|-|y|+\la+\ga+n]$, and viewed as functions of $y_1$. 
 It is  easy to verify that both $f_L$ and $f_R$ are theta functions of order $n+1$ and norm $2|x|+\ga+n-y_2-\dots-y_n$. That the singularities at $[y_1-y_j]=0$ are removable follows from a symmetry argument, and is also apparent from \eqref{wi}.
Thus, it suffices to verify \eqref{mti} for $n+1$ independent values of $y_1$. 
We might as well consider the 
$2n$ values $y_1=x_i$, $y_1=x_i+1$, $1\leq i\leq n$.  By symmetry, it is enough to take $i=1$. On the left, we apply Lemma \ref{sl}. On the right, we observe that  if $y_1=x_1$, all terms with $1\in S$ vanish, and if $y_1=x_1+1$, all terms with $1\notin S$ vanish. In both cases, the sum can be expressed as a sum over subsets of $\{2,\dots,n\}$. In this way, \eqref{mti} 
is  reduced to an equivalent identity with $n$ replaced by $n-1$. 
\end{proof}

A particularly interesting case of \eqref{wi} is   $\gamma=|y|-|x|$, when it
 takes the form
$$Z_n(x;y;\lambda)=\sum_{k=1}^n C_k\frac{[\la+|y|-|x|+n-k]}{[\la+n-k]}, $$
with $C_k$
independent of $\lambda$. Moreover, if
$\eta=1/N$ with $N=2,\dots,n$, we can even write (using that $[x+N]=-[x]$)
$$Z_n(x;y;\lambda)=\sum_{k=1}^{N-1}  D_k\frac{[\la+|y|-|x|+n-k]}{[\la+n-k]}. $$
The existence of such elliptic partial fraction expansions (in the sense of \cite{r})  gives very precise information on $Z_n$ as a function of $\lambda$.

\begin{corollary}\label{lpc}
As a function of $\lambda$,
$$Z_n(x;y;\lambda)\prod_{j=1}^n[\la+j-1] $$
is a theta function of order $n$ and norm $|x|-|y|-\binom n2$. Moreover, if $\eta=1/N$, with  $N=2,\dots,n$,
$$Z_n(x;y;\lambda)\prod_{j=1}^{N-1}[\la+n-j] $$
is a theta function of order $N-1$ and norm $|x|-|y|+n-\binom N2$.
In particular, in any case
$Z_n$ has only single poles in  $\lambda$.
\end{corollary}

These facts are quite remarkable since,   in general, the  individual terms in \eqref{ss} have poles of high multiplicity. The first part of Corollary~\ref{lpc} is also clear from Theorem \ref{wft}. However, in the case $\eta=1/N$, the expression given there still has apparent multiple poles.

\section{Sums of determinants}
\label{sds}

We proceed to show that  in the important case $q^N=1$, the partition function can be expressed as a sum of $N-1$ determinants.

 When dealing with this  case, we
find it 
 convenient to change to multiplicative notation. 
To this end,
we observe that the expression
$$q^{ n(|x|+|y|)/2}Z_n(x;y;\lambda) $$
is invariant under translations by $1/\eta$ in each of the variables $x_i$, $y_i$ and $\lambda$. This follows from Lemma \ref{tfc} and Corollary \ref{lpc}. Thus, there exists a function $\tilde Z_n(x;y;\lambda) $
such that
$$\tilde Z_n(q^{x_1},\dots,q^{x_n};q^{y_1},\dots,q^{y_n};q^\lambda)
=q^{n(|x|+|y |)/2}Z_n(x;y;\lambda).
 $$

We may  write Theorem \ref{wft} as
\begin{multline}\label{mtim}
\tilde Z_n(x;y;\la)=\frac{q^{\frac 12n(n-1)}\prod_{i,j=1}^n x_i\,\theta(y_j/x_i)}{\theta(q)^{n(n-1)}\prod_{j=1}^n\theta(\la q^{j-1})}\\
\times\sum_{\sigma\in S_n}\prod_{1\leq i<j\leq n}\frac{\theta(qy_{\sigma(j)}/y_{\sigma(i)},y_{\sigma(j)}/x_iq)}{\theta(y_{\sigma(j)}/y_{\sigma(i)},y_{\sigma(j)}/x_i)}
\prod_{j=1}^n\frac{\theta(\lambda q^{n-j}y_{\sigma(j)}/x_j)}{\theta(y_{\sigma(j)}/x_j)}.
\end{multline}
Replacing also $q^\ga$ by $\ga$ and recalling the notation \eqref{nn},
Theorem \ref{mt} can similarly be written
\begin{multline}\label{mmi}\tilde Z_n(x;y;\la)=
\frac{(-1)^{\binom n2}\theta(\la q^n)}{\theta(q)^{n^2}\tha(\ga)^n Y\theta(X\la\ga q^n/Y)}\frac{\prod_{i,j=1}^ny_j^2 \theta(x_i/y_j,q x_i/y_j)}{\prod_{1\leq i<j\leq n}x_jy_j\theta(x_i/x_j,y_i/y_j)}\\
\times\sum_{S\subseteq \{1,\dots,n\}}(-1)^{|S|}
\frac{\theta(\la\ga q^{n-|S|})}{\theta(\la q^{n-|S|})}\det_{1\leq i,j\leq n}\left(\frac{\theta(\ga x_i^S/y_j)}{\theta(x_i^S/y_j)}\right),
\end{multline}
where,  in contrast to \eqref{xs},
$$x_i^S=\begin{cases}x_iq,& i\in S,\\ x_i,& i\notin S.\end{cases} $$

We will need the following result, where we use the standard notation
$$(p;p)_\infty=\prod_{j=1}^\infty(1-p^j). $$

\begin{lemma}\label{tdl}
When $N$ is a positive integer, 
$$\frac{\theta(ax;p)}{\theta(x;p)}
=\frac{(p^N;p^N)_\infty^2\,\theta(a;p)}{(p;p)_\infty^2\,\theta(x^N;p^N)}
\sum_{k=0}^{N-1}x^k\frac{\theta(ax^Np^k;p^N)}{\theta(ap^k;p^N)}.
 $$
\end{lemma}

\begin{proof}
Start from the Laurent expansion
\begin{equation}\label{rs}\frac{(p;p)_\infty^2\theta(ax;p)}{\theta(a,x;p)}=\sum_{k=-\infty}
^\infty \frac{x^k}{1-ap^k},\qquad |p|<|x|<1,\end{equation}
which is a special  case of Ramanujan's ${}_1\psi_1$ summation \cite[Eq.\ (II.29)]{gr}.
Replace the summation index $k$ by $k+Nj$, where $0\leq k\leq N-1$ and $j\in\mathbb Z$. Observe that the sum in $j$ can be evaluated using another instance of \eqref{rs}. The restriction $|p|<|x|<1$ is removed by analytic continuation. 
\end{proof}

We now assume that $q^N=1$, with $N=2,3,\dots$. Applying 
Lemma~\ref{tdl} with $a=\ga$, $x=\la q^{n-|S|}$, the sum in \eqref{mmi} takes the form
\begin{multline*}\frac{(p^N;p^N)_\infty^2\,\tha(\ga)}{(p;p)_\infty^2\,\theta(\la^N;p^N)}\sum_{S}(-1)^{|S|}\sum_{k=0}^{N-1}\la^k q^{(n-|S|)k}
\frac{\theta(\ga\la^N p^k;p^N)}{\theta(\ga p^k;p^N)}\\
\times \det_{1\leq i,j\leq n}\left(\frac{\theta(x_i^S\ga/y_j)}{\theta(x_i^S/y_j)}\right)\\
=\frac{(p^N;p^N)_\infty^2\,\tha(\ga)}{(p;p)_\infty^2\,\theta(\la^N;p^N)}\sum_{k=0}^{N-1}\la^k q^{nk}
\frac{\theta(\ga\la^N p^k;p^N)}{\theta(\ga p^k;p^N)}
\\
\times
\det_{1\leq i,j\leq n}\left(\frac{\theta(\ga x_i/y_j)}{\theta(x_i/y_j)}
-q^{-k}\frac{\theta(q\ga x_i/y_j)}{\theta(qx_i/y_j)}
\right).
\end{multline*}

We thus arrive at the following result.

\begin{corollary}\label{fsd} When $q^N=1$,
\begin{multline*}\tilde Z_n(x;y;\la)=\frac{(-1)^{\binom n2}(p^N;p^N)_\infty^2\,\theta(\la q^n)}{(p;p)_\infty^2\,\theta(q)^{n^2}\tha(\ga)^{n-1}\theta(\la^N;p^N)Y\theta(X\la\ga q^n/Y)}\\
\times\frac{\prod_{i,j=1}^ny_j^2 \theta(x_i/y_j,q x_i/y_j)}{\prod_{1\leq i<j\leq n}x_jy_j\theta(x_i/x_j,y_i/y_j)}\\
\times\sum_{k=0}^{N-1} \la^kq^{nk}\frac{\theta(\ga\la^N p^k;p^N)}{\theta(\ga p^k;p^N)}
\det_{1\leq i,j\leq n}\left(\frac{\theta(\ga x_i/y_j)}{\theta(x_i/y_j)}
-q^{-k} \frac{\theta(q\ga x_i /y_j)}{\theta(qx_i/y_j)}
\right),
\end{multline*}
where $\theta(x)=\theta(x;p)$.
\end{corollary}

 In particular, choosing $\gamma=p^{-k}\lambda^{-N}$ for some integer $k$,
one term in the sum vanishes and the partition function is expressed
 as a sum of $N-1$ determinants.

For general $q$, the same method expresses the partition function as an infinite sum of determinants. That is, if we use
 \eqref{rs} to expand the quotient $\theta(\la\ga q^{n-|S|})/\theta(\la q^{n-|S|})$ in \eqref{mmi} we obtain the following identity.

\begin{corollary}\label{rc}
Assume that $|p|<|\la q^{k}|<1$, $0\leq k\leq n$. Then, 
\begin{multline*}
\tilde Z_n(x;y;\la)=\frac{(-1)^{\binom n2}\theta(\la q^n)}{(p;p)_\infty^2\,\theta(q)^{n^2}\tha(\ga)^{n-1} Y\theta(X\la\ga q^n/Y)}\\
\times\frac{\prod_{i,j=1}^ny_j^2 \theta(x_i/y_j,q x_i/y_j)}{\prod_{1\leq i<j\leq n}x_jy_j\theta(x_i/x_j,y_i/y_j)}\\
\times\sum_{k=-\infty}^{\infty} \frac{\la^kq^{nk}}{1-\ga p^k}
\det_{1\leq i,j\leq n}\left(\frac{\theta(\ga x_i/y_j)}{\theta(x_i/y_j)}
-q^{-k} \frac{\theta(q\ga x_i /y_j)}{\theta(qx_i/y_j)}
\right).
\end{multline*}
\end{corollary}

\section{Kuperberg's specialization}
\label{kss}

Kuperberg showed how the partition function for the $6$-vertex model 
can be specialized to the generating function (or $t$-enumeration)
$$\sum_{\text{states}}t^N, $$
where $N$ denotes  the
 number of entries equal to $-1$ in the corresponding alternating sign matrix. 
Using the Izergin--Korepin determinant formula, he  computed the $t$-enumeration for $t=1,2,3$. We intend to generalize this for $t=1$ and $t=2$. To simplify the exposition, we first consider the case of general $t$.

\begin{lemma}\label{ksl}
Let, for each state, $N$ denote the
 number of entries equal to $-1$ in the corresponding alternating sign matrix, that is, the number of blocks of the form $\left(\begin{smallmatrix}a&a-1\\a-1&a\end{smallmatrix}\right)$. Then,
\begin{multline*}\tilde Z_n(q^{-1/2},\dots,q^{-1/2};1,\dots,1;\la)\\
=
q^{-\frac 12\binom{n+1}2}t^{-\binom n2}
\sum_{\text{\emph{states}}}t^N\prod_{\text{\emph{blocks}}}\frac{\tha(\la q^{(3a+3b-c-d)/4})}{\tha(\la q^a)},\end{multline*}
where the product is over all blocks $\left(\begin{smallmatrix}a&b\\c&d\end{smallmatrix}\right)$ and where
$$t=q^{-1/2}\frac{\tha(q)^2}{\tha(q^{1/2})^2} .$$
 In particular, for $p=0$,
\begin{multline*}\tilde Z_n(q^{-1/2},\dots,q^{-1/2};1,\dots,1;\la)\\
=q^{-\frac 12\binom{n+1}2}t^{-\binom n2}\sum_{\text{\emph{states}}}t^N\prod_{\text{\emph{blocks}}}\frac{1-\la q^{(3a+3b-c-d)/4}}{1-\la q^a},
 \end{multline*}
where
$$t=q^{1/2}+q^{-1/2}+2. $$
\end{lemma} 

\begin{proof}
By definition,
$$\tilde Z_n(q^{-1/2},\dots,q^{-1/2};1,\dots,1;\la)
=q^{-n^2/4}\sum_{\text{{states}}}\prod_{\text{{blocks}}}
\tilde 
R^{b-a,d-b}_{d-c,c-a}(\la q^a),
$$
with local weights
$$\tilde R^{++}_{++}(\la)=\tilde R^{--}_{--}(\la)=
q^{1/4}\frac{\tha(q^{1/2})}{\tha(q)}
,
 $$
$$\tilde R^{+-}_{+-}(\la)=-q^{-1/4}\frac{\tha(q^{1/2})\tha(\la q)}{\tha(q)\tha(\la)}, \qquad \tilde R^{-+}_{-+}(\la)=-q^{3/4}\frac{\tha(q^{1/2})\tha(\la q^{-1})}{\tha(q)\tha(\la)}, $$
$$\tilde R^{-+}_{+-}(\la)=q^{1/4}\frac{\tha(\la q^{-1/2})}{\tha(\la )}, \qquad \tilde R^{+-}_{-+}(\la)=q^{-1/4}\frac{\tha(\la q^{1/2})}{\tha(\la )}. $$

We  multiply all local weights by $q^{-1/4}\tha(q)/\tha(q^{1/2})$, and accordingly the prefactor by $(q^{1/4}\tha(q^{1/2})/\tha(q))^{n^2}$. 
 We then multiply $\tilde R^{+-}_{+-}$ by $-q^{1/2}$ and  
 $\tilde R^{-+}_{-+}$ by $-q^{-1/2}$. Since these two types of blocks are equinumerous \cite{k}, this does not change the partition function.
Using that there are $N$ factors of type  $\tilde R^{-+}_{+-}$ and 
$n+N$ factors of type   $\tilde R^{+-}_{-+}$, we obtain
\begin{multline*}\tilde Z_n(q^{-1/2},\dots,q^{-1/2};1,\dots,1;\la)\\
=q^{-n/2}\left(\frac{\tha(q^{1/2})}{\tha(q)}\right)^{n^2-n}
\sum_{\text{{states}}}\left(q^{-1/2}\frac{\tha(q)^2}{\tha(q^{1/2})^2}\right)^N
\prod_{\text{{blocks}}}
\hat
R^{b-a,d-b}_{d-c,c-a}(\lambda q^a),
\end{multline*}
where 
$$\hat R^{++}_{++}(\la)=\hat R^{--}_{--}(\la)=1
,
 $$
$$\hat R^{+-}_{+-}(\la)=\frac{\tha(\la q)}{\tha(\la)}, \qquad \hat R^{-+}_{-+}(\la)=\frac{\tha(\la q^{-1})}{\tha(\la)}, $$
$$\hat R^{-+}_{+-}(\la)=
\frac{\tha(\la q^{-1/2})}{\tha(\la )}, \qquad \hat R^{+-}_{-+}(\la)=\frac{\tha(\la q^{1/2})}{\tha(\la )}. $$
Finally, one checks that in each case
$$\hat
R^{b-a,d-b}_{d-c,c-a}(\lambda q^a,x)=\frac{\tha(\la q^{(3a+3b-c-d)/4})}{\tha(\la q^a)}.
$$
\end{proof}

\section{Dynamical enumeration}
\label{es}
 Kuperberg's proof of the alternating sign matrix theorem 
is based on the six-vertex model with $q$ a primitive cubic root of unity. 
It seems interesting to consider the analogous specialization of the 8VSOS model. In the special case
 $p=0$, we have the following result.

\begin{theorem}\label{det}
When $p=0$ and $q=\om=e^{2\pi i/3}$,
\begin{multline}\label{tsp}\tilde Z_n(\om,\dots,\om;1,\dots,1;\lambda)\\
=
\frac{\om^{\binom{n+1}2}}{(1-\la\om^{n+1})(1-\la\om^{n+2})}\left(
 A_n(1+\om^{n}\la^2 )+(-1)^{n} C_n\om^{2n}\la
\right), \end{multline}
where
$$A_n=\prod_{j=1}^n\frac{(3j-2)!}{(n+j-1)!}$$
is the number of alternating sign matrices of size $n$, and
$$C_n=\prod_{j=1}^n\frac{(3j-1)(3j-3)!}{(n+j-1)!}.$$
\end{theorem}

\begin{remark} $C_n$  is the number of cyclically symmetric plane partitions that fit into a cube of size $n$. This was conjectured by Macdonald \cite{m} and proved by Andrews \cite{a}. Moreover, the product $A_nC_n$ equals  the number of 
$2n\times 2n$
half-turn symmetric alternating sign matrices. This was conjectured by Robbins \cite{ro} and proved by Kuperberg \cite{kus}. 
\end{remark}

\begin{remark} If  
$p\neq 0$ and $q=\om$, it follows from Corollary \ref{lpc} and generalities on theta functions that
\begin{multline*}\tilde Z_n(\om,\dots,\om;1,\dots,1;\lambda)\\
=\frac{1}{\tha(\la\om^{n+1},\la\om^{n+2})}\left(
X_n(p)\theta(-\om^n\la^2;p^2)+Y_n(p)\la\,\theta(-p\om^n\la^2;p^2)
\right),
\end{multline*}
with $X_n$ and $Y_n$ independent of $\lambda$. 
We have not yet been able to find 
simple expressions for these functions, which seem to be  natural elliptic analogues of  $A_n$ and $C_n$. As will be explained below, a solution to this problem would compute the partition function for the three-colour model with domain wall boundary conditions.
\end{remark}

\begin{proof}[Proof of Theorem \ref{det}]
Let $p=0$ and $N=3$ in  Corollary \ref{fsd}. Moreover,  let $\gamma=0$, and replace each $x_i$ by $\om x_i$.  
The resulting identity can be written
\begin{multline*}\tilde Z_n(\om x;y;\lambda)=\frac{(-1)^{\binom n2}}{\om^{\binom n2}(1-\om)^{n^2}(1-\la\om^{n+1})(1-\la\om^{n+2})}
\\
\times\frac{\prod_{1\leq i,j\leq n}(y_j-\om x_i)(y_j-\om^2 x_i)}{\prod_{1\leq i<j\leq n}(x_j-x_i)(y_j-y_i)}\left(D_0+\la\om^n D_1+\la^2\om^{2n}D_2\right),
 \end{multline*}
where
$$D_0=\det_{1\leq i,j\leq n}\left(\frac{1}{y_j-\om x_i}-\frac{1}{y_j-\om^2 x_i}\right)
=(\om-\om^2)^n X\det_{1\leq i,j\leq n}\left(\frac {y_j-x_i}{y_j^3-x_i^3}\right)
, $$
$$D_1=\det_{1\leq i,j\leq n}\left(\frac{1}{y_j-\om x_i}-\frac{\om^2}{y_j-\om^2 x_i}\right)=(1-\om^2)^n\det_{1\leq i,j\leq n}\left(\frac {y_j^2-x_i^2}{y_j^3-x_i^3}\right), $$
$$D_2=\det_{1\leq i,j\leq n}\left(\frac{1}{y_j-\om x_i}-\frac{\om}{y_j-\om^2 x_i}\right)
=(1-\om)^n Y\det_{1\leq i,j\leq n}\left(\frac {y_j-x_i}{y_j^3-x_i^3} \right)
. $$

It remains to let $x_i,\,y_i\rightarrow 1$. This was done by Kuperberg for $D_0$ (which is the same as $D_2$), while $D_1$ can be treated by the same method.
Indeed, it follows from \cite[Theorem 16]{k} and also from \cite[Lemma 13]{l} that
\begin{multline*}\lim_{\substack{x_1,\dots,x_n\rightarrow 1\\ y_1,\dots,y_n\rightarrow 1 }}
\frac 1{\prod_{1\leq i<j\leq n}(x_j-x_i)(y_j-y_i)}\det_{1\leq i,j\leq n}\left(\frac{y_j^k-x_i^k}{y_j^l-x_i^l}\right)\\
=\frac{(-1)^{\binom n2}\prod_{i,j=1}^n(k+l(j-i))}{l^n\prod_{i,j=1}^n(n+j-i)}.
\end{multline*}
It is easy to verify that the right-hand side equals 
$3^{-\binom{n+1}2}A_n$ when $(k,l)=(1,3)$, and 
$3^{-\binom{n+1}2}C_n$ when $(k,l)=(2,3)$. 
After simplification, using also that $(1-\om)^2=-3\om$, we arrive at the desired identity.
\end{proof}

Consider now the case 
$p=0$ and $\eta=-2/3$ of Lemma 
\ref{ksl}, so that $q^{-1/2}=q=\om$ and $t=1$.
We find that  the left-hand side of \eqref{tsp} equals
$$\om^{\binom{n+1}2}\sum_{\text{states}}\prod_{\text{blocks}}\frac{1-\la\om^{-c-d}}{1-\la\om^a}. $$
Note that since $c\not\equiv d\mod 3$, the three numbers $-c-d$, $c$ and $d$ are all noncongruent mod 3, so
$$1-\la\om^{-c-d}=\frac{1-\la^3}{(1-\la\om^c)(1-\la\om^d)}.$$
Thus, we have
\begin{multline*}\sum_{\text{{states}}}
\prod_{\text{blocks}}\frac{1}{(1-\la\om^a)(1-\la\om^c)(1-\la\om^d)}\\
=\frac{(1-\la\om^n)(A_n(1+\om^{n}\la^2 )+(-1)^{n} C_n\om^{2n}\la)
}{(1-\la^3)^{n^2+1}}
.\end{multline*}
Since each matrix entry $a$ in the bulk of the height matrix gives rise to three factors $1/(1-\la\om^a)$,  the left-hand side equals
$$\sum_{\text{{states}}}
\prod_{\text{entries}}\frac{1}{(1-\la\om^a)^3}$$
times the correction
$$(1-\la)^4(1-\la\om)^6\dotsm(1-\la\om^{n-1})^6(1-\la\om^n)^5=
\frac{(1-\la^3)^{2n+2}(1-\la\om^n)}{(1-\la\om^2)^2(1-\la\om^{n+1})^2} $$
arising from boundary entries. Thus, we arrive at the following result, which clearly reduces to the alternating sign matrix theorem when $\lambda=0$.

\begin{corollary}\label{ec}
For each state and for $i=0,1,2$, let $k_i$ denote the number of entries
of the height matrix congruent to $i$ modulo $3$. Then,
\begin{multline}\label{tc}\sum_{\text{{\emph{states}}}}
\prod_{i=0}^2\frac{1}{(1-\la\om^i)^{3k_i}}\\
=\frac{(1-\la\om^2)^2(1-\la\om^{n+1})^2(A_n(1+\om^{n}\la^2 )+(-1)^{n} C_n\om^{2n}\la)}{(1-\la^3)^{n^2+2n+3}}.
\end{multline}
\end{corollary}

This result is best understood in terms of $3$-colourings.
It was observed by Lenard \cite{li} that reducing each  entry in the height matrix modulo $3$ gives a bijection from states to colourings of the
$(n+1)\times(n+1)$ square lattice
 with three colours, such that no adjacent squares have the same colour, and such that (in our case) the  boundary condition 
arising from \eqref{dwbc} is satisfied. It is then natural to introduce the partition function
$$\sum_{\text{{states}}} x_0^{k_0}x_1^{k_1}x_2^{k_2}.$$
In the case of periodic boundary conditions, the thermodynamical limit $n\rightarrow\infty$ of the corresponding partition function was computed by Baxter \cite{b0}.
Corollary \ref{ec}  evaluates the function for fixed $n$ in the case of domain wall boundary conditions, when $x_i$ are
 constrained by the relation
\begin{equation}\label{ac}\left(\frac1{x_0}+\frac1{x_1}+\frac1{x_2}\right)^3=\frac{27}{x_0x_1x_2}. \end{equation}
Indeed, this surface is parametrized by $x_i=t/(1-\la\om^i)^3$, the dependence of $t$ being trivial since $k_0+k_1+k_2=(n+1)^2$. The case of general $p$ similarly corresponds to the unconstrained generating function. 

In spite of the constraint \eqref{ac}, Corollary \ref{ec} contains much information about three-colourings. As an example, we  pick out the coefficient of $\lambda$ on both sides of \eqref{tc} to obtain
$$\sum_{\text{{states}}}(3k_0+3k_1\om+3k_2\om^2)=-2A_n(\om^2+\om^{n+1})+(-1)^nC_n\om^{2n}. $$
Let $$K_i=\sum_{\text{{states}}}k_i. $$
Since, for real values of $a$, $b$, $c$,
$$a+b\om+c\om^2=0\quad\Longleftrightarrow\quad a=b=c, $$
we can deduce that
$$3K_0-(-1)^n C_n=3K_1+2A_n=3K_2+2A_n,\qquad n\equiv 0\mod 3, $$
$$3K_0=3K_1=3K_2+4A_n-(-1)^nC_n,\qquad n\equiv 1\mod 3, $$
$$3K_0+2A_n=3K_1-(-1)^nC_n=3K_2+2A_n,\qquad n\equiv 2\mod 3. $$
Together with the relation $K_0+K_1+K_2=(n+1)^2A_n$, this allows us to solve for
 $K_i$. In Corollary \ref{rcc}, we state the result in terms of the probabilities
$$p_i=\frac{K_i}{(n+1)^2A_n}.$$
In other words, $p_i$ denotes the probability that  a random square from a random 3-colouring has colour $i$. We have also used that
$$\frac{C_n}{A_n}=\frac{2\cdot 5\cdots(3n-1)}{1\cdot 4\cdots(3n-2)}
=\frac{(2/3)_n}{(1/3)_n}. $$

\begin{corollary} \label{rcc}
If $n\equiv 0\mod 3$, the probabilities $p_i$ are given by
$$p_0=\frac13+\frac{4}{9(n+1)^2}+(-1)^n\frac{2\cdot(2/3)_n}{9(n+1)^2(1/3)_n}, $$
$$p_1=p_2=\frac13-\frac{2}{9(n+1)^2}+(-1)^{n+1}\frac{(2/3)_n}{9(n+1)^2(1/3)_n}; $$
if $n\equiv 1\mod 3$,
$$p_0=p_1=\frac13+\frac{4}{9(n+1)^2}+(-1)^{n+1}\frac{(2/3)_n}{9(n+1)^2(1/3)_n}, $$
$$p_2=\frac13-\frac{8}{9(n+1)^2}+(-1)^n\frac{2\cdot(2/3)_n}{9(n+1)^2(1/3)_n}; $$
and if  $n\equiv 2\mod 3$,
$$p_0=p_2=\frac13-\frac{2}{9(n+1)^2}+(-1)^{n+1}\frac{(2/3)_n}{9(n+1)^2(1/3)_n}, $$
$$p_1=\frac13+\frac{4}{9(n+1)^2}+(-1)^n\frac{2\cdot(2/3)_n}{9(n+1)^2(1/3)_n}. $$
\end{corollary}

We point out that the asymptotics of $p_i$ in the thermodynamical limit $n\rightarrow\infty$
 are easily investigated using Stirling's approximation. 
In the first approximation, writing
$$\frac{(2/3)_n}{(1/3)_n}=\frac{\Gamma(1/3)}{\Gamma(2/3)}\,n^{1/3}
+O(n^{-5/3})$$
(note that the term of order $n^{-2/3}$ vanishes),
we  obtain in each case 
$$p_i=\frac 13 +(-1)^nC_1n^{-5/3}+C_2 n^{-2}+O(n^{-3}) \qquad (n\rightarrow\infty),$$
with explicit constants $C_1$ and $C_2$.

\section{Dynamical $2$-enumeration}
\label{tes}

When $q=-1$, the partition function  factors explicitly. Indeed, both 
 Theorem~\ref{wft} and
Corollary \ref{fsd} reduce in this case to known elliptic determinant evaluations. This can also be deduced from recent results of Foda et al.\ \cite{fwz} (as remarked in \cite[\S 2.6]{fwz} a special case of the Felderhof-type model considered there corresponds to the free fermion point, i.e.\ $q=-1$, of Baxter's model).

\begin{proposition}\label{sqp}
When $q=-1$,
\begin{multline*}\tilde Z_n(x;y;\lambda)=\frac{1}{2^{n(n-1)}}\frac{(p;p)_\infty^{2n(n-1)}}{(p^2;p^2)_\infty^{2n(n-1)}}\frac{X\,\theta((-1)^{n+1}\la Y/X)}{\theta((-1)^{n+1}\la)}\\
\times\prod_{1\leq i<j\leq n}x_iy_i\,\theta(-x_j/x_i,-y_j/y_i).
 \end{multline*}
\end{proposition}

\begin{proof}[First proof]
We start from the case $q=-1$ of \eqref{mtim}. Since
$$
\prod_{1\leq i<j\leq n}\frac{\theta(-y_{\sigma(j)}/y_{\sigma(i)})}{\theta(y_{\sigma(j)}/y_{\sigma(i)})}=\operatorname{sgn}(\sigma)\prod_{1\leq i<j\leq n}\frac{\theta(-y_{j}/y_{i})}{\theta(y_{j}/y_{i})},
$$
it can be written
\begin{multline*}\tilde Z_n(x;y;\lambda)=\frac{(-1)^{\binom n2}X^n}{\theta(-1)^{n(n-1)}
\prod_{j=1}^n\theta((-1)^{j-1}\lambda)}\prod_{1\leq i<j\leq n}\frac{\theta(-y_{j}/y_{i})}{\theta(y_{j}/y_{i})}\\
\times\det_{1\leq i,j\leq n}\left(
\prod_{k=1}^{i-1}\tha(-y_i/x_k)\theta((-1)^{n-j}\lambda y_i/x_j)\prod_{k=i+1}^n
\theta(y_i/x_k)
\right).\end{multline*}
By a  generalization of \eqref{fd} due to Tarasov and Varchenko \cite{tv}, 
given more explicitly as  \cite[Corollary 4.5]{rs} (see also
 \cite[Lemma 1]{h}), the determinant equals
\begin{multline*}(-1)^{\binom n2}X^{1-n}\theta((-1)^{n+1}\la Y/X)\prod_{j=1}^{n-1}\theta((-1)^{j-1}\lambda)\\
\times\prod_{1\leq i<j\leq n}x_iy_i\,\theta(-x_j/x_i,y_j/y_i).\end{multline*}
Noting that
$$\tha(-1)=2\frac{(p^2;p^2)_\infty^2}{(p;p)_\infty^2}$$
completes the proof. 
\end{proof}

\begin{proof}[Second proof] Let $N=2$ in 
Corollary \ref{fsd}, and choose $\gamma=p/\la^2$ so that the sum reduces to the term with $k=0$. The determinant is then
$$\det_{1\leq i,j\leq n}\left(\frac{\theta(\ga x_i/y_j)}{\theta(x_i/y_j)}
- \frac{\theta(-\ga x_i /y_j)}{\theta(-x_i/y_j)}
\right).$$
Applying the known identity
$$\frac{\theta(ax)}{\theta(x)}-\frac{\theta(-ax)}{\theta(-x)}
=\frac{2x(p^2;p^2)_\infty^2\,\theta(a,apx^2;p^2)}{(p;p)_\infty^2\,\theta(x^2;p^2)}$$
(which, for instance, follows from \eqref{rs}),
the determinant is evaluated by \eqref{fd}.
\end{proof}

Combining Lemma \ref{ksl} and Proposition \ref{sqp} gives
the following ``dynamical 2-enumeration". When $\lambda=0$, it reduces to 
$$\sum_{\text{{states}}}2^N
=2^{\binom n2}, $$
a result found already in \cite{mrr}.

\begin{corollary}\label{dtc}
Let, for each state, $N$ denote the number of  blocks of the form $\left(\begin{smallmatrix}a&a-1\\a-1&a\end{smallmatrix}\right)$. For $i=0,1,2,3$, let $m_i$ denote the number of blocks  $\left(\begin{smallmatrix}a&b\\c&d\end{smallmatrix}\right)$ with $a+b+c+d\equiv 2i\pmod 8$. 
Then,
\begin{multline}\label{dtci}
\sum_{\text{\emph{states}}}2^N
(1+\la)^{m_0}(1+i\la)^{m_1}(1-\la)^{m_2}(1-i\la)^{m_3}\\
=2^{\binom n2}\cdot\begin{cases}
(1-\lambda^2)^{n^2/2},& n\equiv 0\pmod 4,\\
(1+i\la)
 (1-\lambda^2)^{(n^2-1)/2}
,& n\equiv 1\pmod 4,\\
(1-\lambda)^{(n^2+2)/2}(1+\la)^{(n^2-2)/2},& n\equiv 2\pmod 4,\\
(1-i\la)
 (1-\lambda^2)^{(n^2-1)/2}
,& n \equiv 3\pmod 4.
\end{cases}
\end{multline}
\end{corollary}

\begin{proof}
Let  $q^{-1/2}=i$ in the second part of  Lemma \ref{ksl}, and evaluate the partition function using  Proposition \ref{sqp}. This gives
$$\sum_{\text{states}}2^N\prod_{\text{blocks}}\frac{1-i^{(c+d-3a-3b)/2}\la}{1-(-1)^a\la}=2^{\binom n2}\frac{1+i^{n}\la}{1+(-1)^n\la}
 $$
Since odd and even matrix entries interlace,
$$\prod_{\text{blocks}}(1-(-1)^a\la)=\begin{cases}(1-\lambda)^{n^2/2}(1+\la)^{n^2/2},& n \text{ even},\\
 (1-\lambda)^{(n^2+1)/2}(1+\la)^{(n^2-1)/2},& n \text{ odd},
\end{cases}
 $$
independently of the state. Moreover, since $a+b$ is odd,
$$1-i^{(c+d-3a-3b)/2}\la=1+i^{(a+b+c+d)/2}\la.$$
This yields the desired expression.
\end{proof}

As an indication of the combinatorial meaning of Corollary \ref{dtc}, we identify the coefficient of  $\lambda$ on both sides of \eqref{dtci}. We obtain 
$$
\sum_{\text{{states}}}2^N(m_0+im_1-m_2-im_3)=2^{\binom n2}\cdot\begin{cases}
0,& n\equiv 0\pmod 4,\\
i
,& n\equiv 1\pmod 4,\\
-2,& n\equiv 2\pmod 4,\\
-i
,& n \equiv 3\pmod 4,
\end{cases}
$$
or equivalently
$$
\sum_{\text{{states}}}2^N(m_2-m_0)=\begin{cases}
2\cdot 2^{\binom n2},& n\equiv 2\pmod 4,\\
0
,& \text{else},
\end{cases}
$$
$$
\sum_{\text{{states}}}2^N(m_3-m_1)=\begin{cases}
(-1)^{(n+1)/2}2^{\binom n2}
,& n \text{ odd},\\
0
,& n \text{ even}.
\end{cases}
$$

\end{document}